\documentclass[12pt]{article}
\usepackage{amssymb}
\usepackage{amsmath}
\usepackage{amsfonts}
\usepackage[usenames]{color}
\usepackage{colortbl}
\usepackage{graphicx}

\newcommand{\p}{{\sf P}}
\newcommand{\e}{{\sf E}}
\newcommand{\ind}{{\mathbb{I}}}

\newtheorem{Th}{Theorem}
\newtheorem{Lemma}{Lemma}

\begin{document}
\begin{center}
{\LARGE {\bf Modification of the MDR-EFE method\\
\vskip0,3cm
for stratified samples}}
\end{center}
\vskip0.2cm
\begin{center}
{\large {\bf Alexander Bulinski\footnote{Address: Faculty of Mathematics and Mechanics of the Lomonosov Moscow State University,
Moscow 119991, Russia.
%, Steklov Mathematic\-al Institute of Russian Academy of Science (affiliated
%by RSF grant),
E-mail: bulinski@mech.math.msu.su}, Alexey Kozhevin\footnote{Faculty of Mathematics and Mechanics of the Lomonosov Moscow State University. E-mail: kozhevin.alexey@gmail.com}}}
\end{center}
\vskip0.5cm

The MDR-EFE method of performing identification of relevant factors
within a given collection $X_1,\ldots,X_n$ is developed for
stratified samples in the case of binary response variable $Y$. We
establish a criterion of strong consistency of estimates (involving
$K$-cross-validation procedure and  penalty) for a specified
prediction error function. The cost approach is proposed to compare
experiments with random and nonrandom number of observations.
Analytic results and simulations demonstrate advantages of the
method introduced for stratified samples over that employed for
i.i.d. learning sample.


\vskip0.5cm

\section{Introduction}

The research direction combining probability, statistics and machine
learning for analysis of  mathematical problems of feature selection
is vastly represented in literature along with various applications
of this theory. Quite a number of powerful methods were developed
for different models in the course of such investigations. Several
new variable selection procedures have emerged during the last 20
years. One can refer, e.g., to the following books
%\cite{DQV},
\cite{Ahmed} --
%, \cite{Bolon-Canedo},
\cite{Buhlmann}, \cite{JWHT} --
%, \cite{Koch},
\cite{Marin}, \cite{Ritter}.
Note that many exhaustive, stochastic and heuristic methods
to detect epistasis (in genetics) are considered in \cite{Bul}, \cite{Moore} and \cite{Shang}.

In the paper by M.Ritchie et al. \cite{Ritchie} the {\it multifactor
dimensionality reduction} (MDR) method was proposed to identify the
relevant (in a sense) factors  having influence on a binary response
variable. The review \cite{Gola} demonstrates great popularity of
this method. Some 800 papers published between 2001 and 2014 were
devoted to extensions, modifications and applications of the general
idea suggested in \cite{Ritchie}. The goal of the present paper is
to extend the dimensionality reduction of factors to stratified
samples framework. Here we generalize the approach developed in
\cite{Bulin}--\cite{BulRak} for i.i.d. observations.

Recall some notation. Let $X_1,\ldots,X_n$ be a collection of random
features (or factors) and $Y$ be a response variable depending on
$X:=(X_1,\ldots,X_n)$. We assume that all random elements under
consideration are defined on a probability space
$(\Omega,\mathcal{F},\p)$. Suppose that $X_i$ takes values in a
finite set $\mathbb{X}_i$, $i=1,\ldots,n$. Thus $X\colon\Omega\to
\mathbb{X}$ where $\mathbb{X}:=\mathbb{X}_1\times\ldots\times
\mathbb{X}_n$. We consider $Y$ with values in $\{-1,1\}$ (one uses
also the set $\{0,1\}$).  There are important models in medicine and
biology where $Y$ characterizes the health state of a patient (for
example $Y=1$ means the {\it case}, i.e. a patient is sick, and
$Y=-1$ corresponds to the {\it control}, that is a person is
healthy) and $X$ comprises both genetic and nongenetic factors. The
challenging problem is to predict the risk of certain complex
disease on account of the data $X$ and also to identify the
collection $(X_{i_1},\ldots,X_{i_r})$ of relevant factors ($r<n$)
which are responsible for the disease provoking.

For any $f\colon\mathbb{X}\to \{-1,1\}$  the quality of the forecast of $Y$ by means
of $f(X)$  can be expressed by the following {\it Error function} (see \cite{Bulin}, \cite{BBSSY})
\begin{equation*}\label{d1}
    Err(f) = \mathsf{E}|Y-f(X)|\psi(Y)
\end{equation*}
where a
    {\it penalty function} $\psi\colon \{-1,1\} \to \mathbb{R}_{+}$ is introduced to weight the importance of incorrect prediction of different values of $Y$. The trivial cases $\psi \equiv 0$,
$Y\equiv 0$ or $Y\equiv 1$ are excluded. Clearly,
\begin{equation}\label{e4}
Err(f)=2\sum_{y\in \{-1,1\}}\psi(y)\p(Y=y, f(X)\neq y).
\end{equation}
All the {\it optimal functions}
$f_{opt}$, i.e. $f$ rendering minimum
to $Err(f)$, were described in \cite{Bulin}. It is convenient to take the optimal function
$f^*(x) = \ind\{A\}(x) - \ind\{\overline{A}\}(x)$, $x\in
\mathbb{X}$, where
\begin{equation}\label{A}
    A = \left\{x \in M \colon \mathsf{P}(Y = 1|X = x) > \frac{\psi(-1)}{\psi(-1)+\psi(1)}\right\},
    \end{equation} $M = \{x \in \mathbb{X} \colon \mathsf{P}(X = x) >0\}$, $\ind\{A\}$ stands
for indicator of a set $A$ and $\overline{A}$ means the complement
of $A$. Note that any optimal function $f_{opt}$ gives the same value
to Error function as $Err(f^*)$. In \cite{Bulin} it is explained why
the choice of
\begin{equation}\label{vel}
\psi(y) = \frac{1}{\mathsf{P}(Y = y)},\;\;
y\in \{-1,1\},
\end{equation}
considered in
%by Velez et al.
\cite{Velez} is natural.

However the joint law of $(X,Y)$ is unknown, and therefore $Err(f)$
is unknown ($f^*$ is unknown as well). So it is reasonable to apply
for inference statistical estimates of $Err(f)$ involving
independent identically distributed (i.i.d) random vectors
$(X^1,Y^1),\ldots, (X^N,Y^N),$ having the same law as $(X,Y)$. In
\cite{Bulin}--\cite{BulRak} it is shown how one can use such
estimates to identify the collection of relevant factors. We call
such method MDR-EFE ({\it multifactor dimensionality reduction -
error function estimation}). Moreover, in \cite{B}, \cite{BulRa} and
\cite{BulRak} the asymptotic normality of introduced statistics was
established and in two latter papers a nonbinary response was
studied.

Now we concentrate on the following problem. Assume that the
probability (of disease) $\p(Y=1)$ is rather small. Then in the
sample $(X^1,Y^1),\ldots,(X^N,Y^N)$ for $N$ not too large we could
have too small observations amount with positive response variable
value (equal to one) for sound conclusions. We discuss several
scenarios to overcome this difficulty. Namely, we provide
modifications of some previous results using stratification and also
consider the random number of random observations.

The paper is organized as follows. After the brief Introduction
(Section~1) in Section~2 we prove an auxiliary result concerning the
law of observations having a given response value. Section~3
contains the main result. Here we provide a criterion of strong
consistency of statistics which are the estimates of prediction
error (of approximation of a response $Y$ by means of a function of
factors $X_1,\ldots,X_n$) for stratified samples. We employ the
cross-validation technique or, more precisely, the method of
subsampling and averaging. Along with discussion of established
result we consider an example showing its application. Then we
demonstrate how one can apply the main result for feature selection.
Section~4 is devoted to the cost approach to experiments and the
XOR-model (see \cite{Winh}). Here we compare the MDR-EFE method for
i.i.d. observations and the version of this method for stratified
samples. The simulation results (Section 5) show that even for
rather small samples the proposed version of the MDR-EFE method has
visible advantages.

%simulation of stratified samples according
%to the XOR-model (see \cite{Winh}). Here we compare the MDR-EFE
%method for i.i.d. observations and the version of this method for
%stratified samples employing the cost approach to the plan of
%experiments. The simulation results show that even for rather small
%samples the present version of the MDR-EFE method has visible
%advantages.

    \section{Auxiliary result}

    Let $(X,Y), (X^1,Y^1), (X^2,Y^2), \dots$ be a sequence of i.i.d. random vectors defined on
    a probability space  $(\Omega,\mathcal{F},\mathsf{P})$.
For each $\omega \in \Omega$ we consider a sequence
$Y^1(\omega),Y^2(\omega),\ldots$ and pick all the indices
$1\leqslant j_{-1}^1(\omega) < j_{-1}^2(\omega) < \ldots$ for which
$Y^{j_{-1}^k(\omega)}(\omega)=-1$, $k\in \mathbb{N}$.  In a similar
way we will write all observations $Y^i(\omega)$ with values $1$ as
$\{Y^{j_1^m(\omega)}(\omega)\}_{m\in \mathbb{N}}$ where $1\leqslant
j_1^1(\omega)<j_1^2(\omega)<\ldots$. Recall that for the Bernoulli
trials with probability of success $p$ the Negative binomial random
variable $U_{r,p}$ is introduced as the number of successes needed
to get $r$ failures where $r\in \mathbb{N}$ (one writes $U_{r,p}\sim
NB(r,p)$). Thus
$$
\p(U_{r,p}=k)=\binom{k+r-1}{k}p^k(1-p)^r,\;\;k=0,1,\ldots.
$$
If we consider the events $\{Y^i=1\}$ and $\{Y^i=-1\}$ as a success and failure, respectively  (with
probability $p=\p(Y=1)$ of success),
then
$j_{-1}^r$ has the same law as $U_{r,p}+r$. Therefore
\begin{equation}\label{e2}
\p(j_{-1}^r= m)=
\begin{cases}
\binom{m-1}{m-r}p^{m-r}(1-p)^r,\;\;&m=r,r+1,\ldots,\\
0,\;\;&m=1,\ldots,r-1.
\end{cases}
\end{equation}
For $r=1$ we keep only the first line in \eqref{e2} as in this case
$\{1,\ldots,r-1\}=\varnothing$. By similar reasons $j_1^r$ is
distributed as $U_{r,1-p}+r$ where $U_{r,1-p}\sim NB(r,1-p)$. In
other words $j_1^r$ has the same law as  $G^1_p+\ldots+G^r_p$ where
$G^1_p,\ldots,G^r_p$ are independent random variables having the
Geometric law with parameter $p$ (i.e. $\p(G^1_p=k)=p(1-p)^{k-1}$,
$k=1,2,\ldots$). Thus $\e j_1^r= \frac{r}{p}$ and $j_1^r <\infty$
a.s. for any $r\in \mathbb{N}$ (analogously $\e j_{-1}^r=
\frac{r}{1-p}$ for each $r\in \mathbb{N}$). Note also that one can
find different definitions of Negative binomial and Geometric laws,
that is why we provided the explicit formulae.

Set $Z^k := X^{j^k_{1}}$ for each $k \in \mathbb{N}$. Introduce a collection $\mathcal{B}$ of
all subsets of $\mathbb{X}$. We use the following simple result.

    \begin{Lemma} \label{independence}
        For each $m\in \mathbb{N}$, the random variables $Z^1,\dots,Z^m$ are independent and distributed
        as $X$ given $Y=1$ $($we write $X|Y=1)$, i.e., for any $B\in \mathcal{B}$ and $k=1,\ldots,m$,
$$
\p(Z^k\in B)= \p(X\in B|Y=1).
$$
    \end{Lemma}
    {\it Proof}.
        First of all we show that, for all $B_i \in \mathcal{B}$, $i=1,\ldots,m$,
        \begin{equation}
        \mathsf{P}(Z^1 \in B_1, \dots, Z^m \in B_m) = \mathsf{P}(Z^1 \in B_1) \dots \mathsf{P}(Z^m \in B_m).
        \label{indep}
        \end{equation}
         By the  total probability formula
        \begin{gather*}
        \mathsf{P}(Z^1 \in B_1, \dots, Z^m \in B_m) = \mathsf{P}(X^{j^1_{1}} \in B_1, \dots, X^{j^m_{1}} \in B_m) \\
        = \sum_{\substack{(k_1, \dots, k_m) \in \mathbb{N}^m: \\ k_1 < \dots < k_m}} \mathsf{P}(X^{j^1_{1}} \in B_1, \dots, X^{j^m_{1}} \in B_m, j^1_{1} = k_1 , \dots, j^m_{1} = k_m).
        \end{gather*}
        Note that for arbitrary positive integers $k_1 < k_2 < \dots < k_m$
        \begin{gather*}
        \{j^1_{1} = k_1 , \dots, j^m_{1} = k_m\} = \bigcap_{i=1}^m\{Y^{k_i} = 1\} \cap
\bigcap_{r\in T_m}\{Y^r = - 1\},
        \end{gather*}
        where $T_m = T_m(k_1,k_2,\dots,k_m):= \{1,\dots,k_m\} \setminus \{k_1, \dots, k_m\}$. Next,
        due to
        independence of the random vectors $(X^1,Y^1),(X^2,Y^2),\ldots$, one has
        \begin{gather*}
        \mathsf{P}(Z^1 \in B_1, \dots, Z^m \in B_m)  \\
        = \sum_{\substack{(k_1, \dots, k_m) \in \mathbb{N}^m: \\ k_1 < \dots < k_m}}\!\!\!\mathsf{P}\left(\{X^{k_1} \in B_1,  \dots, X^{k_m} \in B_m, Y^{k_1} = 1 , \dots, Y^{k_m} = 1\}\! \cap\! \left\{
\bigcap_{r\in T_m}\{Y^r \!=\! - 1\}\right\}\right)  \\
        = \sum_{\substack{(k_1, \dots, k_m) \in \mathbb{N}^m: \\ k_1 < \dots < k_m}}\mathsf{P}(X^{k_1} \in B_1,Y^{k_1} = 1)  \dots \mathsf{P}(X^{k_m} \in B_m,Y^{k_m} = 1) \mathsf{P} \left(\bigcap_{r\in T_m}\{Y^r = - 1\} \right)  \\
        =\sum_{\substack{(k_1, \dots, k_m) \in \mathbb{N}^m: \\ k_1 < \dots < k_m}}\prod_{q=1}^m\frac{\mathsf{P}(X^{k_q} \in B_q , Y^{k_q} = 1)}{\mathsf{P}(Y^{k_q} = 1)}
\mathsf{P} \left(\bigcap_{i=1}^m\{Y^{k_i} = 1\} \cap
\bigcap_{r\in T_m}\{Y^r = - 1\}\right)
\end{gather*}
\begin{equation}\label{e1}
        = \prod_{q=1}^m \p(X\in B_q|Y=1) \sum_{\substack{(k_1, \dots, k_m) \in \mathbb{N}^m: \\ k_1 < \dots < k_m}}\mathsf{P}(j^1_{1} = k_1 , \dots, j^m_{1} = k_m)
=\prod_{q=1}^m \p(X\in B_q|Y=1).
\end{equation}

        If we take $B_1 = \dots = B_{k-1} = B_{k+1} = \dots = B_m =  \mathbb{X}$ in \eqref{e1} then
        \begin{gather}
        \mathsf{P}(Z^k \in B_k) =  \mathsf{P}(X \in B_k | Y = 1),\;\;k=1,\ldots,m. \label{single}
%=\mathsf{P}(X^{k} \in B_k | Y^{k} = 1),
        \end{gather}
        Hence, in view of (\ref{single}) the random variables $Z^1,\dots,Z^m$ are identically distributed and
        relation (\ref{indep}) holds.
    $\square$
\vskip0.2cm
{\bf Remark 1}. Lemma 1 is also valid  for i.i.d. vectors $(X,Y),(X^1,Y^1),\ldots,(X^N,Y^N)$ when $X$ takes
values in any space ${\sf S}$ endowed
with some $\sigma$-algebra $\mathcal{B}$ ($X$ is measurable w.r.t. $\sigma$-algebras $\mathcal{F}$ and
$\mathcal{B}$).
In the same way one
        can  prove that $X^{j^1_{-1}},X^{j^2_{-1}},\dots$ are independent and identically
        distributed as $X|Y=-1$.
\vskip0.2cm

    \section{The main result and discussion}
Let $(X^1,Y^1),(X^2,Y^2),\ldots$ be i.i.d. observations having the same law as $(X,Y)$.
Assume that we have a possibility to form a sample with given nonrandom numbers of cases
and controls ($N_1$ and $N_{-1}$ respectively). More exactly, we take
$$
\zeta^1_{N_1}:=\{(X^{j_1^1},1),\ldots,(X^{j_1^{N_1}},1)\},\;\; \zeta^{-1}_{N_{-1}}:=\{(X^{j_{-1}^1},-1),\ldots,(X^{j_{-1}^{N_{-1}}},-1)\}
$$
where the random indices $j_1^k$ and $j_{-1}^k$, $k\in \mathbb{N}$,
were introduced in Section~2. Let $N:=N_1+N_{-1}$ be the size of our
{\it stratified sample} $\zeta_N:=\zeta^1_{N_1} \cup
\zeta^{-1}_{N_{-1}}$. In contrast to the i.i.d. sample
$\xi_N:=\{(X^1,Y^1),\ldots,(X^N,Y^N)\}$ taken from the population
with $law(X,Y)$  we have (according to Lemma~1) two subsamples
$\zeta^1_{N_1}$ and $\zeta^{-1}_{N_{-1}}$ with laws $X|Y=1$ and
$X|Y=-1$, respectively. Thus one cannot use the frequency estimates
(e.g., of $\p(Y=1)$ or $\p(X\in B, Y=1)$ where $B \subset
\mathbb{X}$) constructed by means of $\zeta_N$. Further on we assume
that $N_1=\max\{[aN],1\}$ and $N_{-1}=N-N_1$, here the parameter
$a\in (0,1)$, $N\in \mathbb{N}$ and $[\cdot]$ stands for the integer
part of a number. Suppose that there are some estimators
$\widehat{\mathsf{P}}_{N}^y$ of $\mathsf{P}(Y=y)$ such that
    \begin{equation}
    \widehat{\mathsf{P}}_{N}^y \to \mathsf{P}(Y = y)\text{  a.s.,  } N \to \infty,\;\;y\in \{-1,1\}.\label{Prob}
    \end{equation}
For instance we can assume that $\widehat{\mathsf{P}}_{N}^y$ involve the data
$\{(X^k,Y^k), 1\leqslant k \leqslant \max\{j_1^{N_1},j_{-1}^{N_{-1}}\}\}$. In this case the frequency estimates
of the probabilities $\p(Y=1)$ or
$\p(X\in B, Y=1)$, where $B \subset \mathbb{X}$, mentioned above are strongly consistent
since $\max\{j_1^{N_1},j_{-1}^{N_{-1}}\}\to \infty$ a.s. when $N\to \infty$.
In  Section 4 of the paper we discuss the
%problem of
advantages and disadvantages of employing the stratified samples.
Introduce the vector $\widehat{\mathsf{P}}_{N} :=
(\widehat{\mathsf{P}}_{N}^{-1},\widehat{\mathsf{P}}_{N}^1)$.  Recall
that we exclude the trivial cases  $\mathsf{P}(Y=-1)=0$ or
$\mathsf{P}(Y=1)=0$.

    Let $f_{PA}(x,\zeta_N,\widehat{\mathsf{P}}_{N})$ be a function defining  a \textit{prediction algorithm} i.e. a function with values in $\{-1,1\}$ constructed by means of $x \in \mathbb{X}$,  the sample $\zeta_N$ and $\widehat{\mathsf{P}}_{N}$. In fact we consider a family of functions when
instead of $\zeta_N$ we use its subsamples. Thus we write
$f_{PA}(x,\zeta_N(S),\widehat{\mathsf{P}}_{N})$ for $\zeta_N(S):=
\{(X^j,Y^j),j\in S\}$, $S\subset (\{j_1^1,\ldots,j_1^{N_1}\}\cup
\{j_{-1}^1,\ldots,j_{-1}^{N_{-1}}\})$. For any $f\colon\mathbb{X}\to
\{-1,1\}$ we try to find its estimate $f_{PA}$ which is close in a
sense to $f$ and we employ $f_{PA}$ to estimate $Err(f)$. For this
purpose we also apply \textit{K-fold cross-validation} procedure or,
more precisely, subsampling approach.
%\cite{cv}.
At first, for some fixed $K \in \mathbb{N}$ and each $y\in \{-1,1\}$, we consider a partition of the set
$\{j^1_y,\dots, j^{N_y}_y\}$ into $K$ subsets $S_k^y(N_y,\omega)$, $k = 1,\dots,K$, as follows
%: the $k$th subset of the set $\{a_1,\dots,a_r\}$ consists of the elements $a_i$ with
%indexes belonging to
\begin{equation}\label{partition}
    S_k^y(N_y,\omega):=\left\{j^{i}_y(\omega)\colon i\!\in\! \Bigl\{(k\!-\!1)\left[\frac{N_y}{K}\right]\!+\!1, \dots, k\left[\frac{N_y}{K}\right]\mathbb{I}\{k < K\}\! + \!N_y\mathbb{I}\{k\! =\! K\}\Bigr\}\right\}.
\end{equation}
%   A set $\{j^1_{-1},\dots,j^{N_{-1}}_{-1}\}$ is partitioned similarly
%into subsets $S_k^{-1}(N_{-1},\omega)$, $k=1,\ldots,K$.
Finally, we construct $S_k(N,\omega) = S_k^{1}(N_1,\omega)\cup S_k^{-1}(N_{-1},\omega)$ and introduce
    \begin{multline}\label{Err_est}
    \widehat{Err}_K(f_{PA},\zeta_N,\widehat{\p}_N)\!\! :=\!\!  \frac{2}{K}\!\! \sum_{y\in \{-1,1\}}\!\sum_{k=1}^{K}\sum_{j \in S_k^y(N_y)} \!\!\!\! \frac{\widehat{\psi}(y,\zeta_N(\overline{S_k(N)}),\widehat{\mathsf{P}}_{N})
\mathbb{I}\{f_{PA}^j(N,k) \! \ne\! y \} \widehat{\mathsf{P}}_N^y }{\sharp S_k^y(N_y)}
    \end{multline}
    where
$f_{PA}^j(N,k):=f_{PA}(X^j,\zeta_N(\overline{S_k(N)}),\widehat{\mathsf{P}}_N)$ and $\widehat{\psi}(y,\zeta_N(\overline{S_k(N)}),\widehat{\mathsf{P}}_{N})$ is an estimate of $\psi(y)$, $y\in \{-1,1\}$,
constructed by observations $\zeta_N(\overline{S_k(N)})$
%= \bigcup_{y \in \{-1,1\}} \{X^j_y \colon j \in S^y_k(N)\}$
and $\widehat{\mathsf{P}}_N$; $\overline{S_k(N)}= \{j_1^1,\ldots,j_1^{N_1}\}\cup
\{j_{-1}^1,\ldots,j_{-1}^{N_{-1}}\}\setminus S_k(N)$ and $\sharp$ stands for the
cardinality of a finite set.

    We will suppose that, for each $k = 1,\dots,N$,
\begin{gather}\label{a1}
\widehat{\psi}(y,\zeta_N(\overline{S_k(N)}),\widehat{\mathsf{P}}_{N}) \to \psi(y) \text{  a.s.,  } N \to \infty,\;\;y\in \{-1,1\}.
\end{gather}

Introduce $L(x) = \psi(1)\mathsf{P}(X = x, Y = 1) - \psi(-1)\mathsf{P}(X = x, Y = -1)$, $x\in \mathbb{X}$.
    The following result is an analogue of Theorem~1 \cite{Bulin} for stratified samples.

    \begin{Th} \label{theorem}
        Let $\zeta_N$ be a sample described above, $\psi$ be a penalty function,
$f\colon \mathbb{X} \to \{-1,1\}$ be an arbitrary function and $f_{PA}$ define a prediction algorithm.
%Let also $N_1 = [a N]$ for some $a \in (0,1)$.
Assume that \eqref{a1} is valid and there exists a non-empty set $U \subset \mathbb{X}$ such that, for each $x \in U$ and $k = 1,\dots,K$,  relation
        \begin{gather}
        f_{PA}(x,\zeta_N(\overline{S_k(N)}),\widehat{\mathsf{P}}_{N}) \to f(x) \text{  a.s.,  } N \to \infty,\label{conv}
        \end{gather}
        holds.
        Then, for each $a\in (0,1)$ $($with $N_1=\max\{[aN],1\}$, $N_{-1}=N-N_1)$,
        \begin{gather}
        \widehat{Err}_K(f_{PA},\zeta_N,\widehat{\p}_N) \to Err(f) \text{  a.s.,  } N \to \infty,\label{errorconv}
        \end{gather}
        if and only if
\begin{equation}\label{maincond}
\sum_{k=1}^K\sum_{y\in \{-1,1\}} \sum_{x \in \mathbb{X}_y}y\mathbb{I}\{f_{PA}(x,\zeta_N(\overline{S_k(N)}),\widehat{\mathsf{P}}_{N}) = -y \}L(x) \to 0 \text{  a.s.,  }\;\;N\to \infty,
\end{equation}
%       \begin{gather}
%       \sum_{k=1}^K \sum_{x \in \mathbb{X}_+}\mathbb{I}
%\{f_{PA}(x,\zeta_N(\overline{S_k(N)}),\widehat{\mathsf{P}}_{N}) = -1 \}L(x)
%- \notag \\ -\sum_{k=1}^K \sum_{x \in \mathbb{X}_-}\mathbb{I}
%\{f_{PA}(x,\zeta_N(\overline{S_k(N)}),\widehat{\mathsf{P}}_{N}) = 1 \}L(x)
%\to 0 \text{  a.s.,  }\;\;N\to \infty, \label{condition}
%       \end{gather}
         where
\begin{equation}\label{e5}
\mathbb{X}_y = (\mathbb{X}\setminus U)\cap\{x\in \mathbb{X}: f(x) = y\},\;\;y\in \{-1,1\}.
\end{equation}
%
%       \mathbb{X}_+ = (\mathbb{X}\setminus U)\cap\{x\in \mathbb{X}: f(x) = 1\}, \;\;
%       \mathbb{X}_- = (\mathbb{X}\setminus U)\cap\{x\in \mathbb{X}: f(x) = -1\}.

\end{Th}

{\it Proof}.
%We will denote $S_k^y(N_y)$ by $S_k^y(N)$ because $N_y = N_y(N)$.
It suffices to consider relation
        \begin{gather}
        \frac{2}{K}\sum_{k=1}^{K}\sum_{y\in \{-1,1\}}\psi(y)\sum_{j \in S_k^y(N_y)} \frac{\mathbb{I}\{f_{PA}^j(N,k)  \ne y \} \mathsf{P}(Y = y) }{\sharp S_k^y(N_y)} \to Err(f) \text{  a.s.,  } N \to \infty, \label{withoutphi}
        \end{gather}
        since the difference between $\widehat{Err}_K(f_{PA},\zeta_N,\widehat{\p}_N)$ and the left-hand side of \eqref{withoutphi} tends to zero a.s., $N\to \infty$. Namely, we take into account \eqref{Prob},
\eqref{a1} and the inequality
%       \[
%\frac{1}{\sharp S_k^y(N)}      \sum_{j \in S_k^y(N)}
%\mathbb{I}\{f_{PA}(X^j,\zeta_N(\overline{S_k(N)}),
%\widehat{\mathsf{P}}_{N})  \ne y \} \widehat{\mathsf{P}_N^y} \leqslant 1 \text{  a.s.  }
%       \]
%       and
        \[
\frac{1}{\sharp S_k^y(N_y)}     \sum_{j \in S_k^y(N_y)} \mathbb{I}\{f_{PA}^j(N,k)  \ne y \} \leqslant 1
        \]
which is valid for each $y \in \{-1,1\}$, any $k = 1,\dots,K$ and all $N$.
%(clearly, $\widehat{\mathsf{P}_N^y}\leq 1$).

Now we note that in view of \eqref{e4} relation \eqref{errorconv} is equivalent to the following one
%For each
%$y\in \{-1,1\}$ and any $k=1,\ldots,K$
$$
\frac{2}{K}\sum_{k=1}^{K}\!\sum_{y\in \{-1,1\}}\psi(y)\!\left(\sum_{j \in S_k^y(N_y)}\!\!\! \frac{\mathbb{I}\{f_{PA}^j(N,k)  \ne y \} \mathsf{P}(Y = y) }{\sharp S_k^y(N_y)} - \p(f(X)\neq y,Y=y)\right)\to 0\;\;a.s.
$$
when $N\to \infty$.

For any $y \in \{-1,1\}$ and $m \in \mathbb{N}$, set $
\widetilde{X}^{j^m_y} = \mathbb{I}\{f(X^{j^m_y})  \ne y\} -
\p(f(X^{j^m_y})  \ne y). $ Then $|\widetilde{X}^{j^m_y}|\leqslant 2$,
$\e \tilde{X}^{j^m_y} = 0$ for such $y$ and $m$. Fix any $y \in
\{-1,1\}$ and $k\in \{1,\ldots,K\}$. Consider an array
$\mathcal{A}(y,k):=\{\widetilde{X}^{j^m_y}, j^m_y \in S_k^y(N_y),
N_y=N_y(N)\}$ where $N\in \mathbb{N}$. The strong law of large
numbers for arrays (SLLNA) given in  Theorem~2.1 \cite{SLLNM}
applies, e.g., when  $p=2$, $\psi(x)=|x|^{5/2}$, $x\in \mathbb{R}$,
and $k\in \mathbb{N}$ (we keep here the notation $p$, $\psi$ and $k$
used in the mentioned paper \cite{SLLNM} for other objects). Whence
one has
\begin{equation}\label{reasoning}
\frac{1}{\sharp S_k^y(N_y)} \sum_{j \in S_k^y(N_y)}\mathbb{I}\{f(X^j)  \ne y\} \to \mathsf{P}(f(X^{j_y^1})  \ne y)  \text{  a.s.,  } N \to \infty.
\end{equation}
More precisely, in \cite{SLLNM} a triangular array was considered.
The study of our  array $\mathcal{A}(y,k)$ can be reduced to
analysis of a collection of several triangular arrays. Indeed, in
view of \eqref{partition} the contribution to asymptotic behavior of
$\frac{1}{\sharp S_K^y(N_y)}    \sum_{j \in
S_K^y(N_y)}\mathbb{I}\{f(X^j)  \ne y\}$ of the last summands
appearing in $S_K^y(N_y)$ with indices greater than $[N_y/K]K$ is
negligible. Therefore, for each $k $ belonging to $\{1,\ldots,K\}$
we have an array $\widetilde{\mathcal{A}}(y,k)$ with rows containing
$r_1,r_2,\ldots$ elements ($(r_i)_{i\in \mathbb{N}}$ depends on
$(N_y(N))_{N\in \mathbb{N}}$) such that $r_1\leqslant r_2\leqslant
\ldots$ ($\widetilde{\mathcal{A}}(y,k)= \mathcal{A}(y,k)$ for
$k=1,\ldots,K-1$). Moreover, if $r_m=\ldots =r_q$ then $q-m
\leqslant I$ where  $I=[Ka]$. Now we can consider separately each of
$I$  arrays (taking one row among the $I$ rows of equal length) with
strictly increasing numbers of elements in rows. Then we introduce
auxiliary rows, if necessary, with elements $\mathsf{P}(f(X^{j_y^1})
\ne y)$ to get a triangular array and apply the mentioned SLLNA.
Lemma~1 shows that $\mathsf{P}(f(X^{j_y^1})\neq y)= \mathsf{P}(f(X)
\ne y|Y = y)$. Hence
        \[
        \frac{2}{K}\sum_{k=1}^K \sum_{y\in \{-1,1\}}\psi(y)\sum_{j \in S_k^y(N_y)} \frac{\mathbb{I}\{f(X^j)  \ne y \} \mathsf{P}(Y = y) }{\sharp S_k^y(N_y)} \to Err(f) \text{  a.s.,  } N \to \infty.
        \]
Therefore, relation \eqref{withoutphi} is valid if and only if
$$
\sum_{k=1}^{K}\!\sum_{y\in \{-1,1\}}\psi(y)\frac{1}{\sharp S_k^y(N_y)}\!\sum_{j \in S_k^y(N_y)}\!\!\!
Z^j_{N,k}(y) \to 0\;\;\mbox{a.s.},\;\;N\to \infty,
$$
where
$$
Z^j_{N,k}(y)= (\ind\{f_{PA}(X^j,\zeta_N(\overline{S_k(N)}),\widehat{\mathsf{P}}_{N})  \ne y \}  - \mathbb{I}\{f(X^j)  \ne y \}) \mathsf{P}(Y = y), \;\;y\in \{-1,1\}.
$$
Set $F_{N,k}(x,y) := \mathbb{I}\{f_{PA}(x,\zeta_N(\overline{S_k(N)}),\widehat{\mathsf{P}}_{N} )  \ne y \} - \mathbb{I}\{f(x) \ne y\}$. For each $y \in \{-1,1\}$, $N \in \mathbb{N}$ and $k = 1,\dots,K$, we introduce the random variables
        \[
        Q_{N,k,y} = \frac{1}{\sharp S^y_k(N_y)}\sum_{j \in S^y_k(N_y)}F_{N,k}(X^j,y).
        \]
        Then (\ref{withoutphi}) is  equivalent to the following relation
        \[
        \sum_{k=1}^{K} \sum_{y \in \{-1,1\}} \psi(y) \mathsf{P}(Y = y) Q_{N,k,y} \to 0 \text{  a.s.,  } N \to \infty.
        \]
        Note that $Q_{N,k,y} = Q^{(1)}_{N,k,y} + Q^{(2)}_{N,k,y}$ where
        \begin{gather*}
        Q^{(1)}_{N,k,y} = \frac{1}{\sharp S^y_k(N_y)}\sum_{j \in S^y_k(N_y)}\mathbb{I}\{X^j \in U\}F_{N,k}(X^j,y),\\
        Q^{(2)}_{N,k,y} = \frac{1}{\sharp S^y_k(N_y)}\sum_{j \in S^y_k(N_y)}\mathbb{I}\{X^j \notin U\}F_{N,k}(X^j,y).
        \end{gather*}
We have
        \begin{gather*}
        |Q^{(1)}_{N,k,y}| \leqslant \sum_{x \in U}\frac{1}{\sharp S^y_k(N_y)}\sum_{j \in S^y_k(N_y)}|\mathbb{I}\{f_{PA}(x,\zeta_N(\overline{S_k(N)}),\widehat{\mathsf{P}}_{N})  \ne y \} - \mathbb{I}\{f(x) \ne y\}|.
        \end{gather*}
Condition (\ref{conv}) entails
        \begin{gather*}
        \sum_{k=1}^{K} \sum_{y \in \{-1,1\}} \psi(y) \mathsf{P}(Y = y) Q^{(1)}_{N,k,y} \to 0 \text{  a.s.,  } N \to \infty.
        \end{gather*}
%       for each $x \in U$, $k=1,\dots, K$.
        If $U = \mathbb{X}$, then $Q^{(2)}_{N,k,y} = 0$ for all $N,k$,$y$ and (\ref{withoutphi}) holds. Let
        $U \ne \mathbb{X}$. Then
        \begin{gather*}
        V_{N,k} := \sum_{y \in \{-1,1\}}\psi(y)\mathsf{P}(Y = y)Q^{(2)}_{N,k,y} =
\sum_{x \in \mathbb{X}_{-1} \cup \mathbb{X}_1}\sum_{y \in \{-1,1\}}
G_{N,k,y}(x)
        \end{gather*}
where $N\in \mathbb{N}$, $k=1,\dots, K$, the sets $\mathbb{X}_y$ were introduced in \eqref{e5} and
%= (\mathbb{X}\setminus U)\cap\{x\in \mathbb{X}: f(x) = 1\}$,  =
%(\mathbb{X}\setminus U)\cap\{x\in \mathbb{X}: f(x) = -1\}$ and
        \begin{gather*}
        G_{N,k,y}(x) \!:=\! \frac{\psi(y)\mathsf{P}(Y \!=\! y)}{\sharp S^y_k(N_y)}\sum_{j \in S^y_k(N_y)}\!\!\mathbb{I}\{X^j \!=\! x\}\Big(\mathbb{I}\{f_{PA}(x,\zeta_N(\overline{S_k(N)}),\widehat{\mathsf{P}}_{N})  \ne y \} - \mathbb{I}\{f(x) \ne y\}\Big).
        \end{gather*}
Let $R^y_j(x) := \mathbb{I}\{X^j = x\}\psi(y)\mathsf{P}(Y = y)$, $y\in \{-1,1\}$, $j\in \mathbb{N}$. The
following equalities are valid
        \begin{gather*}
        \sum_{x \in \mathbb{X}_1}\sum_{y \in \{-1,1\}}G_{N,k,y}(x)\! =\! \sum_{x \in \mathbb{X}_1}\mathbb{I}\{f_{PA}(x,\zeta_N(\overline{S_k(N)}),\widehat{\mathsf{P}}_{N})\! =\! -1 \}\sum_{y \in \{-1,1\}}\frac{y}{{\sharp S_k^y(N_y)}}\sum_{j \in S^y_k(N_y)}  R^y_j(x),
        \end{gather*}
%where $R^y_j(x) = \mathbb{I}\{X^j = x\}\psi(y)\mathsf{P}(Y^j = y)$,
        \begin{gather*}
        \sum_{x \in \mathbb{X}_{-1}}\sum_{y \in \{-1,1\}}G_{N,k,y}(x) \!=\! -\!\!\sum_{x \in \mathbb{X}_{-1}}\mathbb{I}\{f_{PA}(x,\zeta_N(\overline{S_k(N)}),\widehat{\mathsf{P}}_{N})\! =\! 1 \}\!\sum_{y \in \{-1,1\}}\frac{y}{{\sharp S_k^y(N_y)}}\!\sum_{j \in S^y_k(N_y)}  R^y_j(x).
        \end{gather*}
        Similarly to (\ref{reasoning}) it can be shown that, for each $y \in \{-1,1\}$,  $k = 1,\dots,K$ and $x \in \mathbb{X}$, relation
        \begin{gather*}
        \frac{1}{{\sharp S_k^y(N_y)}}\sum_{j \in S^y_k(N_y)} R^y_j(x) \to J(x,y) \text{  a.s.,  } N \to \infty,
        \end{gather*}
        holds where $J(x,y) = \psi(y)\mathsf{P}(X = x, Y = y)$.

        Thus, for any $\varepsilon > 0$, $x \in \mathbb{X}$, $y\in \{-1,1\}$ and almost all $\omega \in \Omega$,
        there exists $N_1(\omega,\varepsilon,k,y)$ such that
        \begin{gather}\label{Rcondition}
        \left|\frac{1}{{\sharp S_k^y(N_y)}}\sum_{j \in S^y_k(N_y)} \left(R^y_j(x) - J(x,y)\right)\right| < \varepsilon
        \end{gather}
        for $N > N_1(\omega,\varepsilon,k,y)$. We note that  $L(x) = \sum\limits_{y \in \{-1,1\}}yJ(x,y)$.
        Furthermore,
        \begin{gather*}
        V_{N,k} = \sum_{x \in \mathbb{X}_1}\mathbb{I}\{f_{PA}(x,\zeta_N(\overline{S_k(N)}),\widehat{\mathsf{P}}_N) = -1 \}L(x) \notag\\
        + \sum_{x \in \mathbb{X}_1}\mathbb{I}\{f_{PA}(x,\zeta_N(\overline{S_k(N)}),\widehat{\mathsf{P}}_N) = -1 \}\sum_{y \in \{-1,1\}}\frac{1}{{\sharp S_k^y(N_y)}}\sum_{j \in S^y_k(N_y)} y(R^y_j(x) - J(x,y))  \notag\\
        -\sum_{x \in \mathbb{X}_{-1}}\mathbb{I}\{f_{PA}(x,\zeta_N(\overline{S_k(N)}),\widehat{\mathsf{P}}_N) = 1 \}L(x) \notag\\
        - \sum_{x \in \mathbb{X}_{-1}}\mathbb{I}\{f_{PA}(x,\zeta_N(\overline{S_k(N)}),\widehat{\mathsf{P}}_N) = 1 \}\sum_{y \in \{-1,1\}}\frac{1}{{\sharp S_k^y(N_y)}}\sum_{j \in S^y_k(N_y)} y(R^y_j(x) - J(x,y)).
        \end{gather*}

        Taking into account (\ref{Rcondition}), for any $\varepsilon > 0$, we obtain the inequality
        \begin{gather*}
         \Bigg|\sum_{k=1}^K\Big( \sum_{x \in \mathbb{X}_1}\mathbb{I}\{f_{PA}(x,\zeta_N(\overline{S_k(N)}),\widehat{\mathsf{P}}_N)\! =\! -1 \}\sum_{y \in \{-1,1\}}\frac{1}{{\sharp S_k^y(N_y)}}\sum_{j \in S^y_k(N_y)} y(R^y_j(x) \!-\! J(x,y))  \notag \\  - \sum_{x \in \mathbb{X}_{-1}}\mathbb{I}\{f_{PA}(x,\zeta_N(\overline{S_k(N)}),\widehat{\mathsf{P}}_N)\! =\! 1 \}\sum_{y \in \{-1,1\}}\frac{1}{{\sharp S_k^y(N_y)}}\sum_{j \in S^y_k(N_y)} y(R^y_j(x) \!-\! J(x,y)) \Big) \Bigg| < 2 K \varepsilon \sharp \mathbb{X}
        \end{gather*}
        when $N > N_1(\omega,\varepsilon) := \max\limits_{k = 1,\dots,K,y \in \{-1,1\}}N_1(\omega,\varepsilon,k,y)$.
        Thus, $\sum\limits_{k=1}^{K}V_{N,k} \to 0$ a.s., as $N \to \infty$, if and only if  condition~\eqref{maincond} holds.
$\square$

{\bf Remark 2}. Theorem~1 demonstrates what one has to verify
outside the ``good set'' U (where $f_{PA}$ approximates $f$
point-wise) to guarantee validity of \eqref{errorconv}. To clarify
\eqref{maincond} we can rewrite it (cf. Theorem~1 \cite{BulRak}) in
the following way
$$
\sum_{k=1}^K\sum_{y\in \{-1,1\}} \sum_{x \in \mathbb{X}_y}y\mathbb{I}\{f_{PA}(x,\zeta_N(\overline{S_k(N)}),\widehat{\mathsf{P}}_{N}) \neq f(x) \}L(x) \to 0 \text{  a.s.,  }\;\;N\to \infty. \label{condition}
$$
Note also that condition~\eqref{maincond}
coincides with the corresponding one proposed for the
case of i.i.d. observations without stratification
(cf. Theorem~1 \cite{Bulin}).
%The discussion of \eqref{a1} is provided at the end of this Section.
\vskip0.2cm

    {\bf Remark 3.} According to Corollary 1 \cite{Bulin} the choice
of
\begin{equation}\label{setU}
    U = \left\{x \in M \colon \mathsf{P}(Y = 1|X = x) \ne \frac{\psi(-1)}{\psi(-1)+\psi(1)} \right\},
    \end{equation}
where $M$  appears in \eqref{A}, implies that  (\ref{maincond}) is satisfied and the problem is to prove that
(\ref{conv}) holds on this $U$ for $f$ under consideration and appropriate $f_{PA}$.

    \vskip0.3cm

    {\bf Remark 4.}
     There are various ways to get strongly consistent estimates $\widehat{\p}^y_N$, $y \in \{-1,1\}$. The
     first one is to estimate $\p(Y=y)$, $y\in \{-1,1\}$, by means of another sample. In this case it is
essential to have the samples belonging to the same law.
    The second approach involves  estimates of $\p(Y = y)$ constructed simultaneously with forming of the
    sample
    $\zeta_N$. More exactly, $\widehat{\p}_N^y$ is a frequency estimate of the $\p(Y=y)$ by means of a random
    number of observations $Y^1, Y^2, \dots, Y^{\widetilde{N}}$ where $\widetilde{N} = \max\{j_1^{N_1}, j_{-1}^{N_{-1}}\}$. Such estimate is strongly consistent as $N\to \infty$ since  $\widetilde{N} \geqslant N$ a.s.
We have mentioned in Section~1 that the choice of $\psi$ given by \eqref{vel} is natural.
%   \begin{equation}\label{form_of_psi}
%   \psi(y) = \frac{1}{\mathsf{P}(Y = y)}, \;\; y \in \{-1,1\}.
%   \end{equation}
%Thus for this penalty function we can employ
%   Estimates $\widehat{\psi}(y,\zeta_N(S_k(N)),\widehat{\mathsf{P}}_{N})$
%for such $\psi(y)$ can be constructed on the basis of $\widehat{\mathsf{P}}_{N}$ as
%   $\widehat{\psi}(y,\zeta_N(\overline{S_k(N)}),\widehat{\mathsf{P}}_{N}) := 1/\widehat{\mathsf{P}}_N^y$,
%$y\in \{-1,1\}$.
% (in our construction $\widehat{\mathsf{P}}_y^N\neq 0$).
Moreover, for such $\psi$ we can simplify \eqref{Err_est}. Namely, now it need not contain
$\widehat{\psi}(y,\zeta_N(\overline{S_k(N)}),\widehat{\mathsf{P}}_{N})$ and $\widehat{\mathsf{P}}_N^y$.
%   In virtue of properties of the almost sure convergence they are strongly consistent.
%In the general case, penalty function $\psi$ can be estimated on the basis of
%$\widehat{\mathsf{P}}_{N}$ and estimates for $\p(X \in B|Y = y)$, $y \in \{-1,1\},

\vskip0.3cm
The estimates of  $\p(X \in B, Y = y)$ where $B \subset
\mathbb{X}$ and $y \in \{-1,1\}$ are required if $\psi(y)$ depends
on the joint distribution of $X$ and $Y$. For this estimation we
cannot employ the sample $\zeta_N$. However the Bayes formula can
help to construct the desired estimates. Such approach in the
framework of classification problems was employed by J.Park in
\cite{Park}.

{\bf Example}. We show how for the stratified sample $\zeta_N$ and an optimal function
$f^*$ one can find $f_{PA}$ to employ the result described by Theorem~1.

For $B \subset \mathbb{X}$ the Bayes theorem yields
$$
        \p(Y = 1|X \in B) = \frac{\p(X \in B|Y = 1) \p(Y = 1)}{\p(X \in B|Y = 1) \p(Y = 1)+\p(X \in B|Y = -1) \p(Y = -1)}.
$$
Therefore, one can construct the estimates of $\p(Y = 1|X \in B)$ by means of $\zeta_N$ and $\widehat{\mathsf{P}}_{N}$.
    Plug-in principle suggests us how to get the appropriate $f_{PA}$ for $f^*$. According to this principle we construct estimates of $\mathsf{P}(Y = 1|X = x)$ and $\frac{\psi(-1)}{\psi(-1)+\psi(1)}$. Then we substitute them into the definition of $A$
introduced by \eqref{A}.

For $y\in \{-1,1\}$, $k=1,\ldots,K$ and $N\in \mathbb{N}$, set
$W^y_k(N_y) := \cup_{m \ne k} S^y_m(N_y)$ and assume that \eqref{a1} holds.
%   \begin{equation}\label{psi2}
%   \widehat{\psi}(y,\zeta_N(\overline{S_k(N)}),
%\widehat{\mathsf{P}}_{N}) \to \psi(y) \text{  a.s.,  } N \to \infty,\;\;y\in \{-1,1\}
%   \end{equation}
Consider the following estimates for $\mathsf{P}(Y = 1|X = x)$ and $\frac{\psi(-1)}{\psi(-1)+\psi(1)}$, respectively
\begin{gather}
        g(x,\zeta_N(\overline{S_k(N)}),\widehat{\mathsf{P}}_N):=  \frac{ \widehat{\mathsf{P}}^{1}_N I_1(x,\zeta_N(\overline{S_k(N)}),\widehat{\mathsf{P}}_N)}{\widehat{\mathsf{P}}^{-1}_N I_{-1}(x,\zeta_N(\overline{S_k(N)}),\widehat{\mathsf{P}}_N) + \widehat{\mathsf{P}}^{1}_N I_{1}(x,\zeta_N(\overline{S_k(N)}),\widehat{\mathsf{P}}_N)} , \notag\\
        h(\zeta_N(\overline{S_k(N)}), \widehat{\mathsf{P}}_N):= \frac{\widehat{\psi}(-1,\zeta_N(\overline{S_k(N)}),\widehat{\mathsf{P}}_{N})}{\widehat{\psi}(-1,\zeta_N(\overline{S_k(N)}),\widehat{\mathsf{P}}_{N}) + \widehat{\psi}(1,\zeta_N(\overline{S_k(N)}),\widehat{\mathsf{P}}_{N})}\label{hh}
\end{gather}
    where
\begin{gather*}
    I_{y}(x,\zeta_N(\overline{S_k(N)}),\widehat{\mathsf{P}}_N) := \frac{1}{\sharp W^y_k(N_y)}\sum_{j \in W_k^y(N_y)}\mathbb{I}(X^j = x), \;\; y \in \{-1,1\},\;\;x\in \mathbb{X}.
    \end{gather*}
As usual we set formally $0/0:=0$.
Let us define
\begin{equation*}
        A_N = \{x \in \mathbb{X} \colon g(x,\zeta_N(\overline{S_k(N)}), \widehat{\mathsf{P}}_N) > h(\zeta_N(\overline{S_k(N)}), \widehat{\mathsf{P}}_N) \}
    \end{equation*}
and
    \begin{equation}\label{fPA}
        f_{PA}(x,\zeta_N(\overline{S_k(N)}),\widehat{\mathsf{P}}_N) := \mathbb{I}\{A_N\}(x) - \mathbb{I}\{\overline{A_N}\}(x).
    \end{equation}
    It can be established similarly to the proof of (\ref{reasoning}) that for $y \in \{-1,1\}$
and $x\in \mathbb{X}$
    \begin{equation*}\label{Est_for_pr}
        I_y(x,\zeta_N(\overline{S_k(N)}),\widehat{\mathsf{P}}_N) \to \mathsf{P}(X = x|Y = y) \text{  a.s.,  } N \to \infty.
    \end{equation*}
By virtue of (\ref{Prob}) and (\ref{a1})  we come to the following relations
    \begin{gather*}
    g(x,\zeta_N(\overline{S_k(N)}),\widehat{\mathsf{P}}_N) \to \p(Y=1|X=x) \text{  a.s.,  } N \to \infty, \notag{g} \\
    h(\zeta_N(\overline{S_k(N)}), \widehat{\mathsf{P}}_N) \to \frac{\psi(-1)}{\psi(-1)+\psi(1)}  \text{  a.s.,  } N \to \infty. \label{h}
    \end{gather*}
    For any $x \in A$ one has $\mathsf{P}(Y = 1|X = x) > \frac{\psi(-1)}{\psi(-1)+\psi(1)}$. Therefore, for
    every $\varepsilon>0$ such  that $$\mathsf{P}(Y = 1|X = x) - \varepsilon > \frac{\psi(-1)}{\psi(-1)+\psi(1)} + \varepsilon,$$
    for almost all $\omega \in \Omega$ and  any $x \in A$ there exist $N_2(\omega)$ and $N_3(\omega)$ such
    that
$$
g(x,\zeta_N(\overline{S_k(N)}),\widehat{\mathsf{P}}_N) > \p(Y=1|X=x) - \varepsilon
$$
for each $N > N_2(\omega)$  and
$$
h(\zeta_N(\overline{S_k(N)}),\widehat{\mathsf{P}}_N) < \frac{\psi(-1)}{\psi(-1)+\psi(1)} + \varepsilon
$$
for each $N \!>\! N_3(\omega)$. Thus, for almost all $\omega \!\in\! \Omega$, any $x \!\in\! A$ and
$N \!>\! \max\{\!N_2(\omega),N_3(\omega)\!\}$ the equalities  $f(x) = 1$ and
$f_{PA}(x,\zeta_N(\overline{S_k(N)}),\widehat{\mathsf{P}}_N) = 1$   are true. Similarly, for
$$
x \in U\cap\overline{A} = \left\{x \in M \colon \mathsf{P}(Y = 1|X = x) < \frac{\psi(-1)}{\psi(-1)+\psi(1)} \right\}
$$
it can be shown that for almost all $\omega \in \Omega$ there exists such $N_4(\omega) \in \mathbb{N}$ that,
for $N > N_4(\omega)$, the equalities $f(x) = -1$ and $f_{PA}(x,\zeta_N(\overline{S_k(N)}),\widehat{\mathsf{P}_N}) = -1$ are satisfied.
    Thus condition~(\ref{conv}) holds for $U$ introduced by formula \eqref{setU}. Hence, the estimate of $Err(f^*)$ appearing in Theorem~1 is strongly consistent if $f_{PA}(x,\zeta_N(\overline{S_k(N)}),\widehat{\mathsf{P}}_N)$ is defined by (\ref{fPA}).
The desired $f_{PA}$ is constructed. $\square$

If $\psi$ is introduced by \eqref{vel} then $U$ and $A$ have the
following form
    \begin{gather*}
    U = \left\{x \in M \colon \mathsf{P}(Y = 1|X = x) \ne \p(Y = 1) \right\},\\
    A = \{x \in M \colon \mathsf{P}(Y = 1|X = x) > \mathsf{P}(Y = 1)\}.
    \end{gather*}
Consequently, instead of \eqref{hh} we can set $h(\zeta_N(\overline{S_k(N)}),\widehat{\mathsf{P}}_N) := \widehat{\mathsf{P}}^{1}_N$ and simplify our Example.

    Now we turn to the situation when a response variable $Y$ depends only on the part of factors.
    Let $\{k_1,\dots,k_r\}$ be a subset of $\{1,\dots,n\}$ where $r < n$.
%and $\mathbb{X}_r = \{0,1,2\}^r$.
We say (cf. \cite{Bolon-Canedo}, Ch. 2) that the factors
$X_{k_1},\ldots,X_{k_r}$ are relevant (or collection of  indices
$k_1,\ldots,k_r$ is relevant) whenever for each $x = (x_1,\dots,x_n)
\in M$ it turns out that
\begin{equation}\label{relev}
    \p(Y=1|X_1 = x_1,\dots,X_n = x_n) = \p(Y=1|X_{k_1} = x_{k_1},\dots,X_{k_r} = x_{k_r}).
\end{equation}

For an arbitrary collection $\{m_1,\dots,m_r\} \subset
\{1,\dots,n\}$ set
    \[
    f^{m_1,\dots,m_r}(x) = \begin{cases}
    1, & \p(Y=1|X_{m_1} = x_{m_1},\dots,X_{m_r} = x_{m_r}) > \frac{\psi(-1)}{\psi(-1)+\psi(1)},\;\; x \in M,\\
    -1, & \text{otherwise}.
    \end{cases}
    \]

    If $\{k_1,\dots,k_r\}$ is a relevant collection then the optimal function $f^*$ has the form  $f^{k_1,\dots,k_r}$.
    Therefore, for any relevant collection $\{k_1,\dots,k_r\} \subset \{1,\dots,n\}$ and an arbitrary subset $\{m_1,\dots,m_r\} \subset \{1,\dots,n\}$, the following
    inequality holds
    \[
        Err(f^{k_1,\dots,k_r}) \leqslant Err (f^{m_1,\dots,m_r}).
    \]
    For any $\{m_1,\dots,m_r\} \subset \{1,\dots,n\}$, $x \in \mathbb{X}$ and a penalty function $\psi$,  consider a prediction algorithm with $f^{m_1,\dots,m_r}_{PA}$ such that
$$
    f^{m_1,\dots,m_r}_{PA}(x,\zeta_N(\overline{S_k(N)}),\widehat{\p}_N) = \begin{cases}
    1, & g^{m_1,\dots,m_r}(x,\zeta_N(\overline{S_k(N)}),\widehat{\p}_N) > h(\zeta_N(\overline{S_k(N)}),\widehat{\p}_N),\\
    -1, & \text{otherwise},
    \end{cases}
$$
where $k=1,\ldots,K$, $h$ is defined by \eqref{hh} and
%   Functions $g^{m_1,\dots,m_r}_{PA}$ and $h$ are defined as
$$
    g^{m_1,\dots,m_r}\!(x,\zeta_N(\overline{S_k(N)}),\widehat{\mathsf{P}}_N)\! =\!  \frac{ \widehat{\mathsf{P}}^{1}_N I^{m_1,\dots,m_r}_1(x,\zeta_N(\overline{S_k(N)}),\widehat{\mathsf{P}}_N)}{\widehat{\mathsf{P}}^{-1}_N I^{m_1,\dots,m_r}_{-1}\!(x,\zeta_N(\overline{S_k(N)}),\widehat{\mathsf{P}}_N) \!+\! \widehat{\mathsf{P}}^{1}_N I^{m_1,\dots,m_r}_{1}\!(x,\zeta_N(\overline{S_k(N)}),\widehat{\mathsf{P}}_N)}.
$$
% , \label{g}\\
%   h(\zeta_N(\overline{S_k(N)}), \widehat{\mathsf{P}}_N) =
%\frac{\widehat{\psi}(-1,\zeta_N(\overline{S_k(N)}),
%\widehat{\mathsf{P}}_{N})}{\widehat{\psi}(-1,\zeta_N(\overline{S_k(N)}),
%\widehat{\mathsf{P}}_{N}) + \widehat{\psi}(1,\zeta_N(\overline{S_k(N)}),
%\widehat{\mathsf{P}}_{N})}, \label{h}
%   \end{gather}
Here, for $y \in \{-1,1\}$ and $x\in \mathbb{X}$,
    \begin{gather*}
    I^{m_1,\dots,m_r}_{y}(x,\zeta_N(\overline{S_k(N)}),\widehat{\mathsf{P}}_N) = \frac{1}{\sharp W^y_k(N_y)}\sum_{j \in W_k^y(N_y)}\mathbb{I}(X^j_{m_1} = x_{m_1},\dots,X^j_{m_r} = x_{m_r}).
    \end{gather*}
Set
    \begin{gather*}
    U := \left\{x \in M_{m_1,\ldots,m_r}
\colon \p(Y = 1|X_{m_1} = x_{m_1},\dots,X_{m_r} = x_{m_r}) \ne \frac{\psi(-1)}{\psi(-1)+\psi(1)}\right\}
    \end{gather*}
where $M_{m_1,\ldots,m_r}=\{x\in \mathbb{X}\colon \p(X_{m_1} = x_{m_1},\dots,X_{m_r} = x_{m_r}) > 0\}$.
    Similarly to \cite{Bulin}, it can be shown that for relevant
    collections
    $\{k_1,\dots,k_r\} \subset \{1,\dots, n\}$, arbitrary
    collections
    $\{m_1,\dots,m_r\} \subset \{1,\dots, n\}$,  any $\varepsilon > 0$ and almost all $\omega \in \Omega$,
    \[
        \widehat{Err}_K(f^{k_1,\dots,k_r}_{PA},\zeta_N,\widehat{\mathsf{P}}_N) \leqslant \widehat{Err}_K(f^{m_1,\dots,m_r}_{PA},\zeta_N,\widehat{\mathsf{P}}_N) + \varepsilon \;\; \text{a.s.}
    \]
    when $N$ is large enough. Thus, it is natural to choose as relevant collection of factors
    $X_{k_1},\ldots,X_{k_r}$ that one which has the minimal  prediction error estimate
    $\widehat{Err}_K(f^{k_1,\dots,k_r}_{PA},\zeta_N,\widehat{\mathsf{P}}_N)$.

\section{Cost approach to experiments and the XOR-model}
We compare two approaches concerning the application of the MDR-EFE method for different sample plans.
    The first one was described in \cite{Bulin} and consists in the employment of nonrandom number of i.i.d.
    observations. The second one is considered in this paper and involves the stratified sample. More exactly,
    stratification means the separation of observations taking into account the values of a response variable
    under consideration.
We will denote these methods as iMDR-EFE and sMDR-EFE, respectively.
It seems that the main disadvantage of the second approach is that
too large number of independent observations is required to form a
stratified sample of $N$ elements with fixed cases to controls ratio
when $\mathsf{P}(Y = 1)$ is very small.
%and that number depends on $\mathsf{P}(Y = 1)$. For example, when $\mathsf{P}(Y = 1)$
%is very small
Moreover, we have to skip a lot of observations with $Y^i = -1$.
However it is worth to emphasize that no  information on predictors
$X^i$ is needed at the stage of forming a stratified sample,
therefore we need not get $X^i$ for skipped observations. This is
essential when the measurement of $Y$ is cheaper than the
measurement of $X$ (since $X$ has components $X_1,\ldots,X_n$). Thus
it is not interesting to compare iMDR-EFE and sMDR-EFE for the equal
sizes of samples and we should take into account additional
observations in sMDR-EFE.

    We propose to compare two approaches in the sense of the total cost of experiment. Assume that there is a
    fixed amount of money $C$ ($C \in \mathbb{N}$) for research. Let each observation $(X^i,Y^i)$ cost~$1$ and
    let the ratio of the price of measuring $Y$  to that of~$X$ be $w \in \mathbb{R}_+$. We
    consider the maximal sizes $s_{ind}(C,w)$ and $s_{str}(C,w)$ of the samples which are available in
    experiments organized to apply  iMDR-EFE and sMDR-EFE, respectively. Let $C_{ind}(N,w)$
    be the total cost of the sample having size $N$ and consisting of independent observations. Let
    $C_{str}(N,w)$ be the total cost of the stratified sample of size $N$. Then
    \begin{gather*}
        s_{ind}(C,w) = \max \{N \in \mathbb{N}\colon C_{ind}(N,w) \leqslant C \},\\
        s_{str}(C,w) = \max \{N \in \mathbb{N}\colon C_{str}(N,w) \leqslant C \}.
    \end{gather*}
Clearly, $C_{ind}(N,w) = N$ and therefore $s_{ind}(C,w) = C$. We note that this value is nonrandom and is
known before the experiment.

    Recall (see Remark 4) that $\widetilde{N} = \widetilde{N}(N,a,\omega) = \max\{j_{-1}^{N_{-1}},j_1^{N_1}\}$.
    Hence $$C_{str}(N,w,\omega) = \frac{1}{w+1}N + \frac{w}{w+1}\widetilde{N}(N,a,\omega)$$
    %C_{str}(N,w) =
    because we have to measure $X$ in $N$ observations (which are included into the sample) and to measure $Y$
    in $\widetilde{N}$ observations. Thus $C_{str}(N,w)$ is a random variable which is unknown before the
    experiment. In general case, $s_{str}(C,w)$ is a random variable as well.  However we can select such $N$
    that $C_{str}(N,w)$ does not exceed  $C$ with probability which is not less than
    $1 - \alpha$, $\alpha \in (0,1)$.
We assume that if the cost of experiment exceeds $C$ with given
small probability then some Reserve Foundation will cover the extra
expenses. Now we want to find $s'_{str}(C,w,\alpha)$ so that
    \begin{equation*}
    s'_{str}(C,w,\alpha) = \max\left\{N \in \mathbb{N} \colon \mathsf{P}(C_{str}(N,w) \leqslant C) \geqslant 1-\alpha\right\}.
    \end{equation*}
Obviously, $\widetilde{N}(N,a,\omega) \geqslant N$ for each $\omega \in \Omega$. Therefore for any
$N \in \mathbb{N}$ and $\omega \in \Omega$
    \begin{equation*}
    C_{str}(N,w,\omega) = N + \frac{w}{w+1}(\tilde{N}(N,a,\omega) - N) \geqslant N.
    \end{equation*}
Consequently, for any $N > C$,
    \begin{gather*}
    \mathsf{P}(C_{str}(N,w) \leqslant C) \leqslant \mathsf{P}(C_{str}(N,w) < N) = 0.
    \end{gather*}
So, for arbitrary $C\in \mathbb{N}$, $\alpha\in (0,1)$ and $w>0$, we
observe that
    \begin{equation}\label{conf}
    s'_{str}(C,w,\alpha) = \max\left\{N \in \{1,2,\dots,C\} \colon \mathsf{P}(C_{str}(N,w) \leqslant C) \geqslant 1-\alpha\right\}.
    \end{equation}
It is worth mentioning that $s'_{str}(C,w,\alpha)$ is a non-random
variable which depends on certain parameters and the distribution of
$\widetilde{N}$. We will apply the following result.

    \begin{Lemma}\label{lemma2}
        For each $m \in \{0,1,2\dots\}$ one has
        \begin{equation}\label{NB}
        \mathsf{P}(\widetilde{N} = m) = \begin{cases}
        0, & m < N, \\
        \mathsf{P}(\eta_{-1} = m - N_{-1}) + \mathsf{P}(\eta_{1} = m - N_{1}), & \text{otherwise}
        \end{cases}
        \end{equation}
where $\eta_{-1} \sim NB(N_{-1},\mathsf{P}(Y = 1))$,    $\eta_{1} \sim NB(N_{1},\mathsf{P}(Y = -1))$.
    \end{Lemma}

    {\it Proof}.
    $\mathsf{P}(\widetilde{N} = m) = 0$ for each $m \in \{0,1,2,\dots,N-1\}$ since $\widetilde{N}(N,a,\omega) \geqslant N$ for each $\omega \in \Omega$. If  $m \in \{N,N+1,\dots\}$ then
    \begin{gather*}
    \mathsf{P}(\widetilde{N} = m) = \mathsf{P}\left(\max\{j_{-1}^{N_{-1}},j_1^{N_1}\} = m\right)  \\
    = \mathsf{P}\left(j_{-1}^{N_{-1}} = m, j_1^{N_1} < m\right)
    + \mathsf{P}\left(j_{-1}^{N_{-1}} < m, j_1^{N_1} = m\right)  \\
    + \mathsf{P}\left(j_{-1}^{N_{-1}} = m, j_1^{N_1} = m\right).
    \end{gather*}

    Let $\widetilde{S}_y(k) = \sharp\{i \in \{1,\dots,k\} \colon Y^i = y\}$. We note that $\widetilde{S}_y(k)+\widetilde{S}_{-y}(k) = k$ for each $k \in \mathbb{N}$. Thus, for each $m \in \{N,N+1,\dots\}$ and any $y\in \{-1,1\}$,
    \begin{gather*}
        \left\{ j_{y}^{N_{y}} = m \right\} = \left\{\widetilde{S}_y(m) = N_y\right\}\cap\left\{\widetilde{S}_y(m-1) = N_y-1\right\}
        \subset   \left\{\widetilde{S}_y(m-1) = N_y-1\right\} \\
        = \left\{\widetilde{S}_{-y}(m-1) = m - N_{y}\right\}
        \subset \left\{\widetilde{S}_{-y}(m-1) \geqslant N_{-y}\right\} \subset \left\{ j_{-y}^{N_{-y}} < m \right\}.
    \end{gather*}
Moreover, $j_{-1}^{N_{-1}}(\omega) \ne j_1^{N_1}(\omega)$ for any $\omega \in \Omega$.
    Therefore
    \begin{gather*}
    \mathsf{P}(\widetilde{N} = m) = \mathsf{P}(j_{-1}^{N_{-1}} = m)
    + \mathsf{P}(j_1^{N_1} = m).
    \end{gather*}
According to Section~2, $j_{y}^{N_{y}}$ has the same distribution as $N_y + \eta_y$. Thus (\ref{NB}) holds.
    $\square$

    {\bf Remark 5.}
    Lemma~\ref{lemma2} shows that the law of $\widetilde{N}$ depends on the known parameters $N_{-1}$, $N_{1}$
    and, in general, unknown $\p(Y=1)$. However, if $\mathsf{P}(Y = 1)$ is known or we have its estimates
    constructed by means either of $\xi_{\widetilde{N}}$ or another sample (see Remark 4) then
    $s'_{str}(C,w,\alpha)$ can be evaluated or estimated.

    \vskip0.2cm

    Using Lemma~\ref{lemma2} we can rewrite $\mathsf{P}\left(C_{str}(N,w) \leqslant C\right)$ as
    \begin{gather}
    \mathsf{P}\left(C_{str}(N,w) \leqslant C\right) = \mathsf{P}\left(\widetilde{N}(N,a,\omega) \leqslant \frac{w+1}{w}C - \frac{N}{w} \right)  \notag\\
    = \p\left( N_{1} \leqslant \eta_{-1} \leqslant \frac{w+1}{w}C - \frac{N}{w} - N_{-1}\right) + \p\left(N_{-1} \leqslant \eta_{1} \leqslant \frac{w+1}{w}C - \frac{N}{w} - N_{1}\right) \label{NB2}
    \end{gather}
    where $\eta_{-1}$ and $\eta_{1}$ are defined above.
Thus, for any $N$, the probability $\mathsf{P}\left(C_{str}(N,w) \leqslant C\right)$    can be evaluated or
estimated (see Remark 5).
    Moreover,
    \begin{equation}\label{monot}
    \mathsf{P}\left(C_{str}(N,w) \leqslant C\right) \geqslant \mathsf{P}\left(C_{str}(N+1,\omega) \leqslant C\right)
    \end{equation}
    for all $N \in \mathbb{N}$ since, for each $\omega \in \Omega$,
    \[
    \frac{1}{w+1}N + \frac{w}{w+1}\tilde{N}(N,a,\omega) \leqslant \frac{1}{w+1}(N+1) + \frac{w}{w+1}\widetilde{N}(N+1,a,\omega).
    \]
Formulae (\ref{NB2}) and (\ref{monot}) suggest the following
algorithm. If $p:=\p(Y=1)$ is known we evaluate
$\mathsf{P}\left(C_{str}(N,w) \leqslant C\right)$ by (\ref{NB2}) for
each $N=1,2,\ldots$ (till $C$) and identify the maximal~$N$
belonging to $\{1,\ldots,C\}$ such that \eqref{conf} holds. In this
way we determine $s'_{str}(C,w,\alpha)$. If $p$ is unknown but there
exists its estimate $\widehat{p}_{N'}$ such that $\widehat{p}_{N'} \to p$
almost surely as $N' \to \infty$ (for instance, it may be an
estimate by means of another sample) then for each fixed $C$, $N$
and $w$ we can estimate $\mathsf{P}\left(C_{str}(N,w) \leqslant
C\right)$ by
$$
 \nu_{N'}(N,w,C)\!\! =\!\!  \sum_{N_{1} \leqslant m \leqslant L_1}\!\! \binom{N_1\!+\!m\!-\!1}{m} (1-\widehat{p}_{N'})^{N_{-1}}\widehat{p}_{N'}^{\;m} + \sum_{N_{-1} \leqslant m \leqslant L_{-1}} \!\! \binom{N_{-1}\!+\!m\!-\!1}{m} \widehat{p}_{N'}^{N_1}(1-\widehat{p}_{N'})^{m}
$$
where $L_y:= \frac{w+1}{w}C - \frac{N}{w} - N_{-y}$, $y\in \{-1,1\}$.
Indeed, $\nu_{N'}(N,w,C)$ tends to the right-hand side of (\ref{NB2}) almost surely as $N' \to \infty$. Thus we
can introduce
    \begin{equation*}
    \widehat{s'}_{str,N'}(C,w,\alpha) := \max\left\{N \in \{1,2,\dots,C\} \colon \nu_{N'}(N,w,C) \geqslant 1-\alpha\right\}
    \end{equation*}
and apply the approach proposed for known $p$.
%If $\mathsf{P}\left(C_{str}(N,w) \leqslant C\right) \ne 1-\alpha$ for each $N$
%then $\widehat{s'}_{str,N'}(C,w,\alpha)(\omega) = s'_{str}(C,w,\alpha)$ for almost all
%$\omega \in \Omega$ if $N'$ is greater than some $N'_1(\omega)$.
%We can implement the algorithm similar to the described above in order
%to compute $\widehat{s'}_{str,N'}(C,w,\alpha)$.

\vskip0.2cm
To compare iMDR-EFE and sMDR-EFE we turn to the popular XOR-model (see, e.g., \cite{Winh}).
This model is used in genetics to describe epistasis without main effects.
%   We compare iMDR-EFE and sMDR-EFE for the same law of $(X,Y)$ and the same
%fixed $C$. Recall that in such case the size of the sample for iMDR-EFE is greater
%than for sMDR-EFE.
Namely, let $\mathbb{X} = \{0,1,2\}^n$ ($0,1,2$ correspond to the
number of minor alleles of a specified gene). Assume now that the
components of a random vector $X = (X_1,\dots,X_n)$ are independent
and, for each $i \in \{1,\dots,n\}$, there exists such $p_i \in
(0,0.5]$ that
\begin{equation}\label{x2}
\mathsf{P}(X_i = 0) = (1-p_i)^2,\;\;\mathsf{P}(X_i = 1) = 2p_i(1-p_i),\;\;\mathsf{P}(X_i = 2) = p_i^2.
\end{equation}
This situation is typical for genome-wide association studies (GWAS)
where each $X_i$ corresponds to a single nucleotide polymorphism (SNP)  and $p_i$ is a minor allele
frequency (MAF). Suppose that collection of relevant factors
(describing a binary response variable $Y$) is
$X_{k_1},\ldots,X_{k_r}$ where $\{k_1,\dots,k_r\} \subset
\{1,\dots,n\}$. One also says that collection of indices
$\overline{k}:=\{k_1,\dots,k_r\}$ is relevant.

    In our simulations we employ the XOR-model of dependence between predictors and response variable. It is
    a generalization of the XOR-model described in \cite{Winh} to the case of more than 2 relevant factors.
    Namely, for each $x = (x_1,\dots,x_n)\in \mathbb{X}$,
\begin{equation}\label{x1}
    \mathsf{P}(Y = 1|X = x) = \begin{cases}
    \gamma, & (x_{k_1} + \ldots + x_{k_r}) \text{ mod } 2 = 1, \\
    0, & \text{otherwise}
    \end{cases}
\end{equation}
where $\gamma \in (0,1)$ and $p_{k_1} = \dots = p_{k_r} = 0.5$. Thus
$(X_{k_1},\ldots,X_{k_r})$ is a collection of relevant factors
according to \eqref{relev}. XOR-model is interesting as it possess
the property described by Lemma~\ref{XOR} below.

For any $\overline{s} := \{s_1,\dots,s_q\} \subset \{1,\dots,n\}$, $q\leqslant n$, set
$X_{\overline{s}}:=(X_{s_1},\ldots,X_{s_q})$.
We write $X_{\overline{s}}=x_{\overline{s}}$
where $x_{\overline{s}} = (x_{s_1},\dots,x_{s_q})\in \{0,1,2\}^q$ if
$X_{s_i}=x_{s_i}$ for all $i=1,\ldots,q$.

    \begin{Lemma}\label{XOR}
        Let $\mathbb{X},X,Y$ and $\overline{k}$ be introduced above $($the dependence between $X$ and $Y$ is
        described by \eqref{x1}$)$. Then, for any collection
        $\overline{m} = \{m_1,\dots,m_l\} \subset \{1,\ldots,n\}$,  a response variable $Y$ is dependent with $X_{\overline{m}}$
if and only if $\overline{k}\subset \overline{m}$.
    \end{Lemma}

{\it Proof}. Formula \eqref{x1} shows that $Y$ and $X_{\overline{m}}$ are dependent
if $\overline{k}\subset \overline{m}$. Let now $\overline{m}=\overline{u}\cup\overline{v}$ where
$\overline{u}\subsetneq \overline{k}$
and $\overline{v}\subsetneq \{1,\ldots,n\}\setminus \overline{k}$. Introduce $t:=\sharp \overline{u}$.
Since
    $
    \p(X_{\overline{m}} = x_{\overline{m}}) \ne 0
    $ for each $x_{\overline{m}} \in \{0,1,2\}^{l}$ it remains to establish that \begin{equation}\label{indepxor}
    \p(Y=1|X_{\overline{m}} = x_{\overline{m}}) = \p(Y=1).
    \end{equation}
Evidently
\begin{gather}
    \p(Y=1|X_{\overline{m}} = x_{\overline{m}})  \\
 = \sum_{x_{\overline{k} \setminus \overline{u}} \in \{0,1,2\}^{r-t}} \p(Y = 1|X_{\overline{u}} = x_{\overline{u}}, X_{\overline{k} \setminus \overline{u}} = x_{\overline{k} \setminus \overline{u}},X_{\overline{v}}=x_{\overline{v}}) \p(X_{\overline{k} \setminus \overline{u}} = x_{\overline{k} \setminus \overline{u}}) \\
    = \gamma \sum_{x_{\overline{k} \setminus \overline{u}} \in \{0,1,2\}^{r-t}} \mathbb{I}\left\{ \left( \sum_{i \in \overline{k} \setminus \overline{u}} x_i + \sum_{i \in \overline{u}} x_i
\right) \text{ mod } 2 = 1  \right\}  \p(X_{\overline{k} \setminus \overline{u}} = x_{\overline{k} \setminus \overline{u}}) \\
    = \gamma \p \left( \sum_{i \in \overline{k} \setminus \overline{u}}
X_i \text{ mod } 2 \ne \sum_{i \in \overline{u}} x_i \text{ mod } 2  \right) \label{x3}
\end{gather}

In view of \eqref{x2} the laws of $(X_1,\ldots,X_n)$ and
$(X_1',\ldots,X_n')$ coincide where
$
    X'_i = \sum_{j=1}^{2}B_i^j
$
and $B_1^1,B_1^2,\ldots,B_n^1,B_n^2$ are independent Bernoulli variables such that
$\p(B_i^k\!=\!1)=\!p_i$, $k=1,2$, $i\!=\!1,\ldots,n$.
Hence, $law(\sum_{i\in \overline{s}}X_i)\!\! =\!\! law(\sum_{i\in \overline{s}}(B_i^1\!+\!B_i^2))$
for any $\overline{s}\!\subset\! \{1,\ldots,n\}$.
Thus, for $\overline{s}:=\{s_1,\ldots,s_q\}\subset \overline{k}$ ($q\leqslant r$), one has
$$
\p \left( \sum_{i \in \overline{s}} X_i \text{ mod } 2 = 0 \right) = \p\left( \sum_{i = 1}^{2q} \tilde{B}_i \text{ mod } 2 = 0  \right) = \sum_{j=0}^{q} C_{2q}^{2j} 2^{-2q} = \frac{1}{2}
$$
where $\widetilde{B}_1,\ldots,\widetilde{B}_{2r}$ are i.i.d. Bernoulli random variables with
probability of success $0.5$. Therefore,
$$
\p \left( \sum_{i \in \overline{s}} X_i \text{ mod } 2 = 1 \right)= \frac{1}{2}.
$$
Taking into account \eqref{x3} we observe that
$\p(Y=1|X_{\overline{m}} = x_{\overline{m}}) = \frac{\gamma}{2}$ for
any $x_{\overline{m}}\in \{0,1,2\}^l$. To complete the proof note
that
    \begin{gather*}
    \p(Y = 1) = \sum_{x \in \mathbb{X}} \p(Y = 1| X = x) \p(X = x) \\ = \sum_{x_{\overline{k}} \in \{0,1,2\}^r} \p(Y = 1| X_{\overline{k}} = x_{\overline{k}}) \p(X_{\overline{k}} = x_{\overline{k}})   = \gamma \p\left(\sum_{i \in \overline{k}} X_i \text{ mod } 2 = 1\right) = \frac{\gamma}{2}.
    \end{gather*}
Thus \eqref{indepxor} holds. $\square$

    \textbf{Remark 6.}
    Lemma 3 shows that only the whole collection of relevant factors $X_{k_1},\ldots,X_{k_r}$
and not some of its subcollections determines the response $Y$ in
XOR-model.
    \vskip0.2cm

To complete this Section we discuss some asymptotical properties of
a random sample as $C\to \infty$. For a given $C$ the iMDR-EFE
method involves $C$ observations whereas its sMDR-EFE counterpart
(using the same amount of money $C$) leads to $N$ observations and
\begin{equation}\label{main}
\frac{1}{w+1}N + \frac{w}{w+1}\max\{j_{-1}^{N-([aN]\vee 1)},j_1^{[aN]\vee 1}\} \leqslant C,
\end{equation}
here $a\in (0,1)$.
Clearly,  \eqref{main} implies that $N\leqslant C$.
Now we consider the problem whether one can take $N=\lambda C$ for some
$\lambda \in (0,1)$ such that with probability close to one
inequality \eqref{main} is satisfied.

\begin{Lemma}
For an arbitrary $($fixed$)$ $a\in (0,1)$ and $w $ inequality \eqref{main} is valid with probability
tending to one as $C\to \infty$ if and only if
\begin{equation}\label{ineq}
\lambda < \lambda_0:= (1+w)\left(1+ w\max\left\{\frac{a}{p},\frac{1-a}{1-p}\right\}\right)^{-1}.
\end{equation}
\end{Lemma}

{\it Proof}. We can rewrite \eqref{main} in the following way
%$$
%\frac{w}{w+1}\max\{j_{-1}^{[\lambda C]-[a[\lambda C]]},j_1^{[a[\lambda C]]}\} \leq C- \frac{[\lambda C]}{w+1},
%$$
%i.e.
\begin{equation*}\label{eq1}
\max\{j_{-1}^{[\lambda C]-[a[\lambda C]]},j_1^{[a[\lambda C]]}\} \leqslant \frac{C(w+1)-[\lambda C]}{w}.
\end{equation*}
Note that
\begin{gather*}
\p\left(j_1^{[a[\lambda C]]} \leqslant \frac{C(w+1)-[\lambda C]}{w}\right)
\\
= \p\left(\frac{j_1^{[a[\lambda C]]}- \frac{[a[\lambda C]]}{p}}{\sigma\sqrt{[a[\lambda C]]}}\leqslant \frac{1}{\sigma\sqrt{[a[\lambda C]]}}
\left(\frac{C(w+1)-[\lambda C]}{w} - \frac{[a[\lambda C]]}{p}\right)\right)
\end{gather*}
where $\sigma^2:= \frac{1-p}{p^2}$ (the variance of a random variable following the Geometric
law with parameter $p$). The central limit theorem for i.i.d. random variables
having finite variance yields that
$j_1^{[a[\lambda C]]} \leqslant \frac{C(w+1)-[\lambda C]}{w}$  with probability
tending to one as $C\to \infty$ if and only if $(w+1-\lambda)/w - a\lambda/p>0$, that is
%$$
%\frac{w+1-\lambda}{w}- \frac{a\lambda }{p} >0,
%$$
%that is
\begin{equation}\label{z1}
\lambda < (1+w)\left(1+\frac{aw}{p}\right)^{-1}.
\end{equation}
In a similar way we can claim that
$$
\p\left(j_{-1}^{[\lambda C]-[a[\lambda C]]} \leqslant \frac{C(w+1)-[\lambda C]}{w}\right)\to 1,\;\;C\to \infty,
$$
if and only if
\begin{equation}\label{z2}
\lambda < (1+w)\left(1+\frac{(1-a)w}{1-p}\right)^{-1}.
\end{equation}
Thus the statement of Lemma~1 is valid if and only if
\eqref{z1} and \eqref{z2} are true simultaneously, i.e. relation \eqref{ineq} holds.
$\square$

Thus the optimal boundary $\lambda_0$ is found.

\section{Simulations}
    We will consider different levels of the total cost $C$ and parameter $\gamma$ of XOR-model.
Note that now $\p(Y=1)=\frac{\gamma}{2}$, i.e. $p= \frac{\gamma}{2}$. We employ  two levels of parameter $w$.
The value $w=0$ corresponds to the extreme case when values of a response variable
are obtained free of charge. This situation can be viewed as the limit one when the price of the measurement
of $Y$ is very low w.r.t. the price of $X$ measurement. According to  \cite{Velez} the application of MDR method
is reasonable for balanced  datasets, i.e. with equal number of cases and controls. Therefore, in the case of
the stratified samples we consider $a=0.5$ and only even $N$ in \eqref{conf}.
For each combination of these parameters we generate $D$ independent datasets for iMDR-EFE and for sMDR-EFE
applications in order to evaluate performance of these methods by Monte Carlo experiments.

    \begin{table}[ht]
        \begin{center}
            \renewcommand{\arraystretch}{1.2}
            \begin{tabular}[b]{|c||c|}
                \hline
                \textbf{Parameter} & \textbf{Value} \\
                \hline
                \hline
                $\psi(y)$ & $\left(\p(Y=y)\right)^{-1}$\\
                \hline
                $D$& 100\\
                \hline
                $n$ & 100 \\
                \hline
                $K$ & 5 \\
                \hline
                $r$ & 3 \\
                \hline
                ($k_1,\dots,k_r$) & (2,3,5) \\
                \hline
                $a$ & $0.5$\\
                \hline
                $C$ & 5 levels: $\{100, 200, 300, 400, 500\}$ \\
                \hline
                $w$ & 2 levels: $\{0, 0.1\}$\\
                \hline
                $\gamma$ & 3 levels: $\{0.05, 0.1, 0.2\}$ \\
                \hline
                $\alpha$ & 0.05 \\
                \hline
            \end{tabular}
            \caption{\label{tab:setting}Simulation scheme parameters.}
        \end{center}
    \end{table}

The values of parameters used in our simulations are given in Table \ref{tab:setting}.
For relevant factors $X_{k_1},\ldots,X_{k_r}$ we set $p_{k_i} = 0.5$, $i=1,\ldots,r$ (according to XOR-model) whereas for other (non-relevant) factors $X_i$ the corresponding $p_i$ are drawn independently from the uniform distribution $U(0.05,0.5)$ to generate each dataset. For our datasets collections of $p_i$ are drawn independently. Such simulation setting for $p_i$ was proposed in \cite{Dehman}. The interval $[0.05,0.5]$ for uniform law was taken since $p_i$ is considered as MAF.

%In order to generate $Y$ given $X$ we generate Bernoulli variable such that success
%probability is equal to $0$ if $X_{k_1}+\dots +X_{k_r} \text{ mod } 2 = 0$ and
%is equal to $\gamma$ if $X_{k_1}+\dots +X_{k_r} \text{ mod } 2 = 1$.

%In order to generate $Y$ given $X$ we generate Bernoulli variable $B(\gamma)$ which
%is independent of X such that success probability is equal to $\gamma$. Then
%we put $Y = -1$ if $X_{k_1}+\dots +X_{k_r} \text{ mod } 2 = 0$  and $Y = 2 B(\gamma) - 1$
%if $X_{k_1}+\dots +X_{k_r} \text{ mod } 2 = 1$.

We take a Bernoulli variable $B(\gamma)$ having the success probability $\gamma$.
Let $B(\gamma)$ and $X$ be independent.
Introduce $Y:= 2 B(\gamma) \ind\{\sum_{v=1}^r X_{k_v}\; \mbox{mod}\; 2 =1\} -1$ where $X_{k_1},\ldots,X_{k_r}$ are the relevant factors. Then one can verify that \eqref{x1}
holds. In such a way we generate independent vectors $(X^j,Y^j)$, $j=1,2,\ldots$.

Taking $w = 0.1$ we will assume that $p$ is known in order to use it
for determining $s'_{str}(C,w,\alpha)$ introduced in \eqref{conf}.
For $w = 0$ the size of the sample $s'_{str}(C,0,\alpha)$ is equal
to $C$ as in the case of independent observations. Then it is
non-random and fixed before the experiment and we estimate $\p(Y=1)$
by means of $Y^1,\ldots,Y^{\widetilde{C}}$. Here $\widetilde{C} = \widetilde{C}(C,a,\omega) = \max\{j_{-1}^{C_{-1}},j_1^{C_1}\}$, $C_1 = [aC]$ and $C_{-1} = C - [aC]$. In order to measure the
method performance power we use $TMR$ (true model rate) which is
defined as
$$
TMR:= \frac{1}{D}\sum_{d=1}^{D}T_d
$$
where
$$
T_d = \begin{cases}
1 ,& \text{if all relevant factors are identified correctly in $d$-th dataset},\\
0, & \text{otherwise.}
\end{cases}
$$

Below the phrase ``implementation of iMDR-EFE (sMDR-EFE) for sample'' means that we select such collection of
$r$ factors which has the minimal $\widehat{Err}_K$.

The simulation procedure can be described in the following way.
%\begin{enumerate}

I. Fix some values $\gamma$ and $C$ from Table \ref{tab:setting}, other parameters except $w$ are fixed as well.

II. Calculate $s'_{str}(C,w,\alpha)$ where $w = 0.1$ and $\alpha = 0.5$ assuming that $p = \mathsf{P}(Y=1)$ is known (and equals $\gamma/2$).

III. Perform $D$ times
    \begin{enumerate}
        \item generation of independent $p_i\sim U(0.05,0.5)$ for non-relevant factors;
        \item generation of independent $X^j$, having law introduced by $\eqref{x2}$, and corresponding $Y^j$ until we have 3 samples:
            \begin{enumerate}
                \item sample $\xi_C$ which consists of $C$ independent observations;
                \item stratified sample $\zeta_{N}$ which consists of $N = s'_{str}(C,w,\alpha)$ observations (for $w = 0.1$);
                \item stratified sample $\zeta_{C}$ which consists of $C$ observations (for $w = 0$). While  generating observations $(X^j,Y^j)$ in order to form $\zeta_{C}$ we estimate $p$ (and write ``$p$ is estimated by means of observations $Y^1,\ldots,Y^{\widetilde{C}}$'');
            \end{enumerate}
        \item implementation of iMDR-EFE for $\xi_C$, $p$ is unknown ($\psi(y)$ as well);
        \item implementation of iMDR-EFE for $\xi_C$, $p$ is known therefore $\psi(y)$ is known and $$\widehat{\psi}(y,\xi_C(\overline{S_k(C)})) = 1/\mathsf{P}(Y=y),$$
            here we use the notation $\widehat{\psi}(y,\xi_C(\overline{S_k(C)}))$ of \cite{Bulin} for i.i.d. observations;
        \item implementation of sMDR-EFE for $\zeta_N$ where $N = s'_{str}(C,w,\alpha)$, $p$ is known therefore $\widehat{\mathsf{P}}_{N} = (1-p,p)$ and $$\widehat{\psi}(y,\zeta_N(\overline{S_k(N)}),\widehat{\mathsf{P}}_{N}) = 1/\mathsf{P}(Y=y);$$
        \item implementation of sMDR-EFE for $\zeta_C$, $p$ is unknown and estimated by means of observations $Y^1,\ldots,Y^{\widetilde{C}}$;
        \item implementation of sMDR-EFE for $\zeta_C$, $p$ is known thus $\widehat{\mathsf{P}}_{C} = (1-p,p)$ and $$\widehat{\psi}(y,\zeta_C(\overline{S_k(C)}),\widehat{\mathsf{P}}_{C}) = 1/\mathsf{P}(Y=y);$$
        \item assignment of $T_d$ to every implementation of the method under consideration.
    \end{enumerate}
IV. Evaluate $TMR$ for each approach.
%\end{enumerate}

\begin{figure}[!ht]
                        \centering
                        \includegraphics[width=1\linewidth]{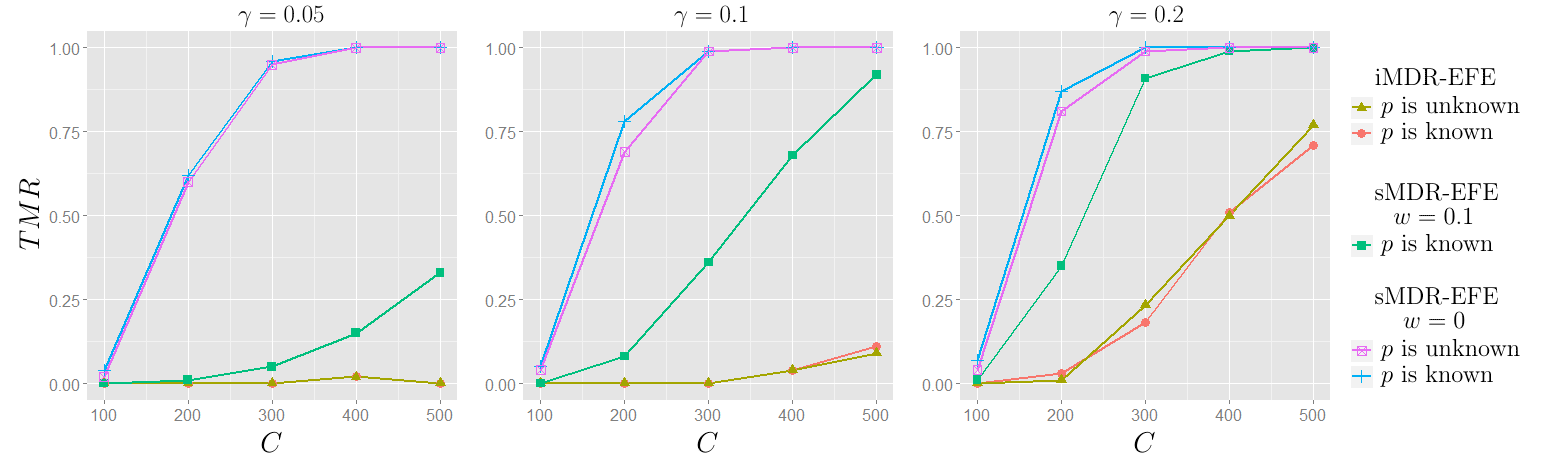}
                        \caption{Performance (TMR) of iMDR-EFE and sMDR-EFE for different levels of $\gamma=2p$.}
                        \label{fig:big_plot1}

                        \centering
                        \includegraphics[width=0.85\linewidth]{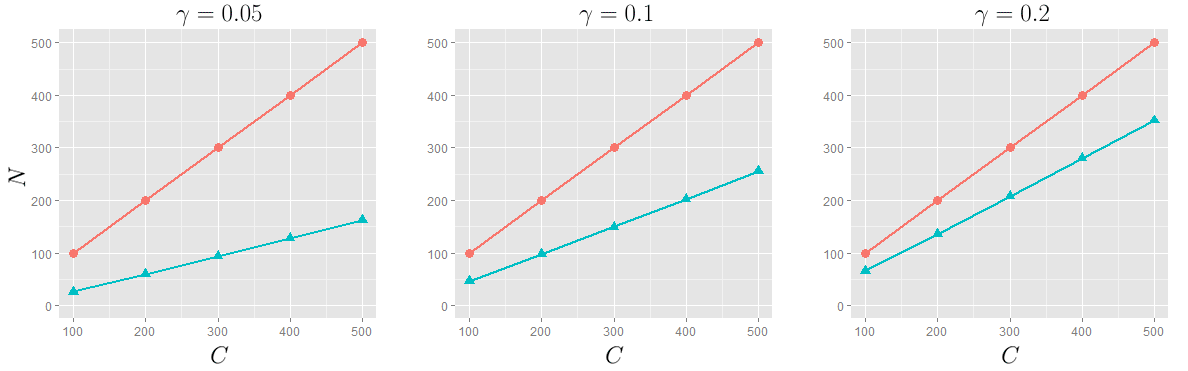}
                        \caption{Sizes of the samples depending on the cost $C$ for different levels of $\gamma=2p$. The red line corresponds to iMDR-EFE and sMDR-EFE with $w = 0$, the blue line corresponds to sMDR-EFE with $w = 0.1$.}
                        \label{fig:big_plot2}
                              \end{figure}

%   The case of $w = 0$ corresponds to the situation of the free measurement of
%response $Y$, thus it is practically impossible. However,
\vskip0.1cm
    The results of our simulations are shown in Figures \ref{fig:big_plot1} and \ref{fig:big_plot2}. Thus for
    the fixed total cost $C$ a stratified sampling gives better results than independent one. Despite of the
    fact that  $s_{ind}(C,\alpha) > s'_{str}(C,w,\alpha)$ (when $w > 0$, see Figure 2), sMDR-EFE demonstrates
    better performance than iMDR-EFE in all 3 models for each value of $C$ and different scenarios concerning $p$. Moreover, small $w$ permits to include more observations into the
stratified sample. Figure~1 suggests that it leads to better
performance of the sMDR-EFE method.

\vskip0.2cm
%According to Lemma 4,
Let $a=0.5$ and $w=0.1$. Then, for $\gamma$  taking values
$0.05$,  $0.1$ and $0.2$ ($p=\gamma/2)$ Lemma 4 yields that the corresponding values of
$\lambda_0$ are equal to $0.366(6)$, $0.55$ and $0.733(3)$,
respectively. For any $a\in (0,1)$, all $p\in (0,1)$ and $w=0$ we
get $\lambda_0=1$. Clearly if $w$ is close to $0$
then $\lambda_0$ is close to $1$.

\vskip0.2cm
We note that there are various possibilities to take into account
the cost of experiments with random and nonrandom number of
observations. This interesting problem will be considered
separately.

\section{Acknowledgements} The work of A.V.Bulinski is supported by the Russian Science Foundation
under grant 14-21-00162 and performed at the Steklov Mathematical
Institute of Russian Academy of Sciences. The problems setting
belongs to A.V.Bulinski, all theoretical results are established
jointly with A.A.Kozhevin. Computer simulations are carried out by
A.A.Kozhevin.

%is carried out with support of Russian Science Foundation (RSF), grant 14-21-00162.

\end{document}